# An Exact Pricing Algorithm for Revenue Maximization under the Logit Demand Function

Moddassir Khan Nayeem [a], Omar Abbaas [a*], Suzan Alaswad [b], Sinan Salman [b]

[a] Department of Mechanical, Aerospace, and Industrial Engineering, The University of Texas at San Antonio, San Antonio, TX 78249, USA

[b] Zayed University, PO Box 144534, Abu Dhabi, United Arab Emirates

## Abstract

Determining the optimal selling price is a challenge in revenue management, especially in markets characterized by nonlinear and price-sensitive demand. While traditional models, such as linear, power, and exponential demand functions, offer analytical convenience, they often fail to capture realistic purchase dynamics, leading to suboptimal pricing. The logit demand function addresses these limitations through its bounded, S-shaped curve, offering a more realistic representation of consumer behavior. Despite its advantages, most existing literature relies on heuristic approaches, such as pricing at the inflection point, which prioritizes maximum price sensitivity but does not guarantee maximum revenue. This study proposes a novel, exact pricing algorithm that analytically derives the revenue-maximizing price under the logit demand function using the Lambert $W$ function. By providing a closed-form solution, the approach eliminates reliance on heuristic iterative methods and corrects the common practice of considering the inflection point price as market price. In fact, we demonstrate that the optimal price is consistently lower than the inflection-point price under reasonable assumptions, leading to lower prices for consumers and higher revenue for sellers. Numerical experiments illustrate the proposed algorithm and examine the changes in the optimality gap as demand function parameters vary. Results indicate that the optimal price is consistently lower than the inflection-point price, with an average 20% price reduction accompanied by a 15% increase in revenue.

## 1. Introduction

Strategic pricing plays a critical role in maximizing revenue, especially when demand is highly sensitive to price. Firms must identify the optimal selling price by understanding the price-demand relationship. In practice, even marginal price changes can yield substantial financial impact. For instance, a study by McKinsey & Company found that one retailer adjusted over 100,000 prices across 150 SKUs and improved its return-on-sales by 3-5% within just three months through disciplined pricing analytics [1]. Similarly, in e-commerce, companies such as Amazon, which adopt dynamic pricing, have reported revenue increases of 20-30% [2]. Setting prices too high suppresses demand, whereas pricing too low may boost volume but erode margins. The optimal price balances these effects to maximize revenue. In today's competitive, algorithm-driven marketplaces, frequent price adjustments underscore the need for models that accurately capture consumer responses to price changes [3]. However, many analytical pricing and inventory models still rely on oversimplified assumptions about how demand reacts to price [4].

In classical inventory and pricing models, demand is often modeled using linear, power, or exponential functions due to their simplicity and analytical tractability [5]. Although Lau and Lau [6] emphasized that even minor variations in the demand curve can result in major shifts in optimal solutions, the existing literature provides little insight or clear justification for selecting one demand function over another [7, 8]. Phillips [5] argued that neither the linear function nor the iso-elastic function can realistically represent the price-demand relationship. For example, the linear model can predict negative demand at higher prices and yields unrealistic infinite elasticity at its extreme points. Power demand functions assume constant price elasticity, failing to capture shifts in consumer behavior across different price levels. Exponential demand curves can predict market sizes exceeding practical limits when prices drop too low. These behaviors of linear, power, and exponential functions are graphically presented in Figure 1. These drawbacks become especially problematic in contemporary pricing environments, where consumer choices are influenced by nonlinear perceptions of value and bounded rationality. Therefore, among various price-dependent demand models, the logit demand function has gained increasing attention due to its ability to reflect realistic consumer behavior, offer global stability across price ranges, and avoid the implausible extremes of traditional models such as linear or power functions [9, 10].

The logit demand function addresses these issues by modeling demand as a reverse S-shaped curve, initially concave, transitioning through an inflection point, and eventually convex. This shape aligns with empirical observations, capturing how consumers exhibit the highest price sensitivity near a central point of the logit curve and less sensitivity at extremely low or high prices. For a given product, most sellers will price it around the middle point (see Figure 2) [9, 11]. In that region, small price changes can attract or deter price-sensitive customers to a given seller. However, on the left-hand side, if a seller captures most of the market size by offering low prices, a small price change (discount or increase) will not significantly alter the demand since the seller is already offering low prices compared to the market. Customers who are not buying despite the low prices may have other reasons not to buy from this seller, such as higher service expectations or negative views of the brand. Similarly, on the right-hand side of the curve, when a seller has high prices for the product, total demand will be lower, however, existing demand is less likely to change significantly with small price adjustments. Customers who buy despite the high prices may be loyal to the seller, prioritize convenience over price, or have other motivations. Importantly, the logit function bounds demand between zero and a finite market size, making it suitable for capturing realistic demand responses across the entire price spectrum [12, 13].

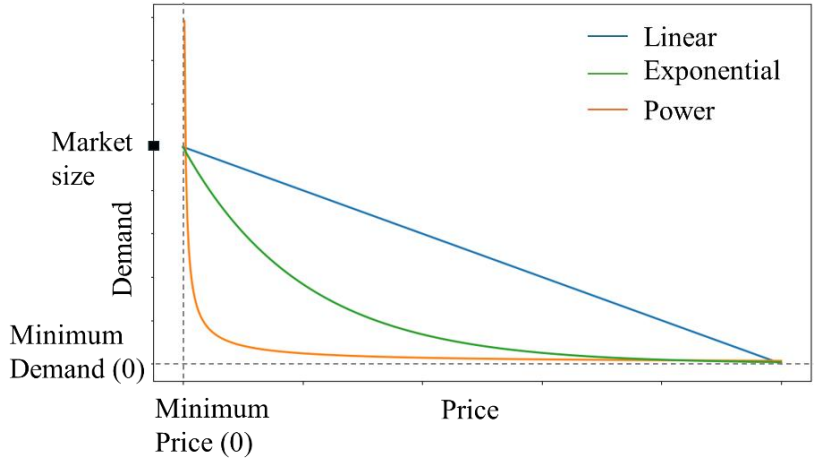

Figure 1: Graphical illustration of different demand functions

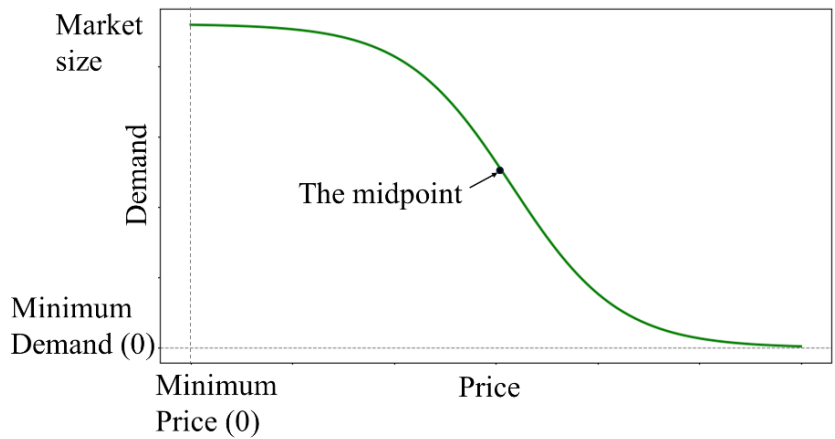

Figure 2: Logit demand function

Despite its advantages, applying the logit function for pricing decisions presents two significant challenges. First, many studies continue to rely on heuristic approaches, such as the inflection point, to guide price-setting decisions. While this point, where demand sensitivity peaks, is often referred to as the "market price" [9, 11], we show in this paper that it does not necessarily align with revenue maximization. The revenue-optimal price may occur earlier on the curve, where the trade-off between unit price and demand volume is more favorable. Second, in practical settings, firms often lack perfect knowledge of the demand function's parameters. Determining the optimal price becomes especially difficult when these parameters must be estimated from sales data. Businesses frequently engage in price experimentation, testing different price points to infer consumer behavior. For brick-and-mortar retailers, this can be accomplished using digital price tags that enable real-time adjustments and data collection [14]. However, many businesses operate under constraints that limit the frequency or range of allowable price changes [15-17]. Too much experimentation risks eroding customer trust, and prolonged periods of suboptimal pricing can result in revenue loss, while too little experimentation yields imprecise parameter estimates. This dilemma gives rise to the well-known learning-earning trade-off, which is similar to the exploration-exploitation dilemma in machine learning and evolutionary optimization [18, 19].

To address these challenges, this study proposes an Exact Pricing Algorithm (EPA) that derives the global revenue-maximizing price under the logit demand function in closed form. Beyond introducing a Lambert $W$ based solution, the paper contributes more broadly by showing that the realistic behavioral demand (through the logit structure) can still yield analytically tractable pricing formulas, bridging the gap between empirical accuracy and mathematical solvability. This stands in contrast to most existing work, which relies on numerical optimization or heuristic benchmarks such as pricing at the inflection point. We further prove that the optimal price always lies strictly to the left of the traditional market price (inflection point), correcting a common practice in the pricing literature and offering a theoretically grounded justification for more competitive price levels. Importantly, this result establishes conditions under which the revenue-maximizing price also aligns with consumer welfare by simultaneously lowering prices and increasing demand. Together, these contributions provide not only a closed-form pricing rule but also meaningful economic insight for real-world pricing decisions in data-driven and competitive markets.

The remainder of this paper is structured as follows. We review the relevant literature in Section 2. Section 3 presents the problem statement, including the properties of the logit demand function and the proposed algorithm for finding the revenue-maximizing price. Section 4 presents the numerical experiments, and Section 5 concludes with key findings and future research directions.

## 2. Literature review

Understanding consumer demand behavior is crucial for developing effective pricing strategies that maximize revenue and enhance market competitiveness. Over the years, various demand models have been proposed to capture the relationship between price and demand. This section first reviews traditional demand and pricing models and their limitations, before introducing more advanced alternatives, particularly the logit demand model.

### 2.1 Traditional demand models

Pricing models form the analytical foundation of revenue management and provide mathematical frameworks for determining optimal prices in response to consumer demand. Early classifications distinguish between static and dynamic models, deterministic and stochastic formulations, and single versus multi-product settings [20, 21]. These models have been widely adopted in industries such as airlines, hospitality, and retail, where price adjustment and capacity control are integrated to maximize revenue [22]. In practice, however, many firms continue to employ traditional pricing approaches such as cost-plus pricing, where a fixed markup is added to production costs [23], value-based pricing, which sets prices according to perceived customer value rather than cost [24], hourly pricing, where charges depend on time spent delivering a service [25], and fixed (flat) pricing, which offers a predetermined total price for a product [26]. While these approaches are widely used for their simplicity and transparency, they generally lack the analytical rigor needed to capture nonlinear consumer responses or optimize revenue dynamically. Consequently, operations and revenue management research has focused on more formal mathematical models, where analytical tractability has traditionally depended on simplified demand functions, most commonly linear, power, or exponential forms. Recent developments in pricing analytics, as discussed by Gallego and Topaloglu [27], have emphasized the integration of choice-based demand models, such as the logit and Multinomial Logit (MNL) formulations, into revenue management frameworks.

Historically, prior to the adoption of choice-based frameworks, the development of pricing strategies in operations research was driven by the search for analytically tractable yet practical demand models. Among the earliest and most influential frameworks is the Economic Order Quantity (EOQ) model, which traditionally assumes a fixed demand rate and decouples pricing decisions from inventory policies [28]. As a result, several demand models have been proposed and extensively studied to bridge this gap, especially within the joint pricing and inventory management literature.

Among traditional models, the linear demand function has been one of the most widely used models in pricing and inventory optimization problems. Its basic form, $d(p) = \alpha - \theta p$, where $\alpha > 0$ is the initial demand value, $\theta > 0$ is the price sensitivity parameter, and $p$ is the price, offers simplicity and ease of parameter estimation. Linear models are straightforward to analyze and can yield closed-form solutions in many optimization settings [29-32]. However, their limitations are well documented. As price increases toward $\alpha/\theta$, demand approaches zero, and beyond this point, it becomes negative, a scenario with no practical interpretation. Furthermore, the linear model implies that price elasticity of demand becomes infinite at the boundary of the feasible price range, resulting in unrealistic elasticity behavior [5, 9, 11].

To address these limitations of the linear model, researchers introduced the power demand function, $d(p) = \theta p^{-l}$ where $\theta > 0$ and $l > 1$ is the constant elasticity. This model assumes constant price elasticity, meaning that demand decreases proportionally at a fixed percentage rate as price increases, and has been widely adopted in pricing and replenishment models due to its analytical tractability [33-36]. While more realistic than linear models in capturing price sensitivity, the power function still suffers from structural issues. Most notably, as the price approaches zero, demand tends to infinity, an implausible scenario in any real-world market [8].

Additionally, its fixed elasticity fails to capture varying consumer responsiveness across price ranges, which is particularly problematic in heterogeneous markets [8, 37].

The exponential demand function, $d(p) = \mu e^{\theta p}$, where $\mu$ is the market size and $\theta$ is price sensitivity, offers an alternative that maintains positivity and smooth decline across price levels. It avoids the negative demand issue of the linear model and the infinite demand problem of the power model. Consequently, the exponential form has been widely adopted in inventory and pricing models to represent either time- or price-dependent demand patterns [38-42]. However, its major drawback lies in its behavior at low price points. As price approaches zero, the exponential function predicts an unrealistic increase in demand, overstating customer sensitivity. In addition, the curve lacks an inflection point and fails to represent varying consumer elasticity across different price regions, making it inadequate for modeling realistic demand patterns observed in empirical data [9].

Table 1 summarizes these demand functions and highlights their mathematical forms, parameter constraints, and qualitative properties. Despite satisfying the basic Law of Demand, where demand decreases as price increases, none of these classical models (linear, power, and exponential) fully capture the nuanced behavior of real-world consumers. For instance, empirical studies in the automotive and food industries have shown that consumers exhibit varying price sensitivities depending on product characteristics, promotional contexts, and psychological thresholds [9, 43-45]. These dynamics are not adequately represented by linear, power, or exponential functions, which are either too rigid or produce implausible results under extreme price conditions. Nonetheless, Duan and Ventura [11] emphasize that many supply chain optimization models continue to rely on these classical forms, often for reasons of analytical tractability rather than empirical accuracy. This trend persisted despite warnings from researchers like Lau and Lau [6], who cautioned that even minor deviations in demand curve assumptions could drastically alter optimal pricing and replenishment policies.

Table 1: Comparison between demand functions

| Demand functions | Formulas | Parameters |
| --- | --- | --- |
| Linear | $d(p) = \mu - \theta p$ | $\mu > 0, \theta > 0$ |
| Power | $d(p) = \theta p^{-l}$ | $\theta > 0, l > 1$ |
| Exponential | $d(p) = \mu e^{-\theta p}$ | $\mu > 0, \theta > 0$ |
| Logit 1 | $d(p) = \left(\frac{\mu e^{-(\theta p)}}{1+e^{-(\theta p)}}\right)$ | $\mu > 0, \theta > 0$ |
| Logit 2 | $d(p) = \mu \left(\frac{e^{-(\alpha+\theta p)}}{1+e^{-(\alpha+\theta p)}}\right)$ | $\mu > 0, \alpha < -2, \theta > 0$ |

## 2.2 Logit demand model

As pricing decisions become more data-driven and behaviorally informed, classical models face significant limitations. Their inability to account for bounded market sizes, dynamic consumer preferences, and realistic elasticity behavior calls for more sophisticated alternatives. This has motivated the use of models like the logit demand function, which offers a bounded, nonlinear, and S-shaped demand curve that aligns more closely with observed market dynamics.

Over the past two decades, the logit demand model has gained prominence as a behaviorally grounded alternative, rooted in discrete choice theory [46-48]. The logit has several appealing properties: it produces a bounded, S-shaped (sigmoid) demand curve, captures variable elasticity, and better aligns with real-world purchase dynamics. Moreover, it models the probability of purchase as a function of price, making demand bounded between zero and the market size. This feature enables the logit model to accommodate diminishing marginal sensitivity and the inflection point where price sensitivity is maximized [8, 9, 11]. Importantly, it avoids predicting negative or infinite demand and can be parameterized to reflect varying consumer segments.

Among logit-based models, the MNL framework is widely known for modeling consumer choice among multiple products. While powerful for assortment optimization and multi-product pricing, the MNL model introduces complexities that are not always necessary in settings focused on a single product or a homogeneous offering [49, 50]. In such cases, the consumer decision can be effectively reduced to a binary outcome: to purchase or not to purchase. This simplification retains the key behavioral advantages of the logit model, such as bounded demand and variable elasticity, while avoiding the additional assumptions and computational burdens of full MNL estimation [47, 51]. Several studies have recognized this simplification as valid and useful for modeling single-product demand when consumer heterogeneity can be addressed through carefully parameterized sensitivity and baseline demand factors [11, 52]. Accordingly, in this study, we adopt a single-product logit demand function that is derived as a special case of the MNL framework, as discussed by Phillips [9]. This approach enables a tractable yet behaviorally realistic representation of demand that is particularly well-suited for revenue optimization. The chosen logit model can accommodate varying price sensitivity and aligns with empirical observations in markets where consumer responses are not uniform across price levels.

Recent literature has applied logit models across various domains, including assortment planning [53-56], dynamic pricing [11, 12, 57-61], and supply chain optimization [12, 13, 61-63]. For instance, Gao et al. [53] and Arhami et al. [56] employed nested and impatient-consumer extensions of the logit model to study product assortment decisions in multi-stage retail settings, demonstrating improved alignment with observed consumer choice patterns. In supply chain contexts, Agrawal and Yadav [61] examined an integrated production-inventory and pricing problem in a single manufacturer-multiple buyers supply chain, incorporating price-dependent demand and proposing profit allocation schemes that balance fairness and practical implementation. Similarly, Ghasemy Yaghin et al. [7] formulated a joint pricing and lot-sizing model for a two-echelon supply chain using a logit-based demand function, showing its capability to capture bounded price sensitivity and inventory interactions more effectively than linear or constant-elasticity models. Suh and Aydin [57] analyzed dynamic pricing strategies for substitutable products with limited inventories using a logit framework that accounts for substitution effects and stock limitations. San-Jose et al. [64] combined a price-logit function with a time-dependent demand pattern to design an inventory model under backlogged shortages. Macias-Lopez et al. [65] included the logit form, among other nonlinear demand types, in a perishable inventory system that jointly accounts for price, stock, and time effects. Similarly, Diaz-Mateus et al. [62] developed a nonlinear optimization framework for a two-echelon supply chain that explicitly considers consumers' maximum willingness to pay, using a constrained logit model to estimate demand across different socioeconomic segments. Complementing these developments, Abbaas and Ventura [12, 13] incorporated logit-based demand structures into procurement auction mechanisms, enhancing bid evaluations and supplier selection under consumer-driven demand variability.

Beyond inventory and pricing contexts, Huang [66] developed a logit-based mode choice model to analyze transport pricing and elastic demand in a bi-modal highway-transit system, while Fu and Wilmot [67] applied a sequential logit formulation to model dynamic travel demand during hurricane evacuation scenarios. Extending the use of logit structures to logistics and infrastructure planning, Bagloee et al. [68] proposed a logit-based mathematical programming framework for facility location in supply chain networks, capturing probabilistic decision-making among suppliers and terminals.

## 2.3 Research gaps

Despite advances in modeling price-dependent demand, a clear research gap remains in providing a robust, closed-form analytical solution for revenue-maximizing pricing under single-product logit demand function. Moreover, the choice of logit specification (between Logit 1 and Logit 2) is also crucial. Prior works by Yaghin et al. [7] and Diaz-Mateus et al. [62] employed Logit 1 formulation (See Table 1) to determine optimal retail prices within two-echelon supply chains, while Fattahi et al. [63] introduced a variant known as Logit 2, utilizing discrete price levels to guide pricing decisions across supply chain networks. A careful examination of the Logit 1 demand curve reveals that it closely resembles the exponential demand function; however, as price approaches zero, the demand only converges to half of the market potential. To address this limitation, this study adopts the Logit 2 demand function proposed by Phillips [9], which offers a more practical and flexible framework for modeling price-dependent demand.

However, even with the behaviorally realistic Logit 2 structure, determining the revenue-maximizing price remains analytically unresolved in much of the existing literature. Most studies continue to rely on heuristic benchmarks, such as the inflection point, which indicates maximum price sensitivity but does not necessarily align with revenue maximization [9, 11]. This reliance highlights the need for a closed-form analytical solution for the Logit 2 model that can directly yield the optimal price without resorting to approximations or iterative procedures.

To fill this gap, the present study introduces a novel pricing algorithm that uses the Lambert $W$ function to derive an exact, closed-form solution for the revenue-maximizing price under the Logit 2 demand model. While Aravindakshan and Ratchford [52] previously applied the Lambert $W$ function to derive equilibrium prices and shares in logit-based utility models, their focus was on market share and strategic interaction, rather than revenue optimization. In contrast, this study is the first to derive a closed-form pricing solution for the Logit 2 model, specifically designed for revenue maximization in a single-product context, offering both theoretical clarity and practical relevance. This contribution provides a rigorous analytical benchmark for data-driven pricing methods and enriches the literature on logit-based demand modeling.

## 3. Problem statement

The logit demand function is widely adopted in the literature to model the price-sensitive demand rate of products, including finished, intermediate, and raw goods, due to its ability to capture nonlinear consumer responses and maintain bounded demand across the full price range. It reflects variable marginal sensitivity to price: sensitivity is low at very low and very high prices and peaks at the inflection point. This pattern aligns with observed purchasing behavior and avoiding the unrealistic extrapolations of linear, power, or exponential forms. The price-sensitive logit demand function is presented in Section 3.1, and the revenue-maximizing algorithm (EPA) is presented in Section 3.2.

### 3.1 Price-sensitive logit demand function

The logit demand function proposed by Phillips [9] is expressed in equation (1):

$$d(p) = \mu \left( \frac{e^{-(\alpha+\theta p)}}{1 + e^{-(\alpha+\theta p)}} \right) \quad (1)$$

where $\mu > 0$ represents the overall market size. The parameter $\theta$ is the price sensitivity coefficient, and it must be greater than zero ($\theta > 0$) for the demand to be decreasing with increasing price $p$. The location parameter $\alpha$ shifts the sigmoidal curve horizontally along the price axis. Greater values of $\alpha$ shift the curve to the left, while smaller values shift it to the right. Practically $\alpha$ should be less than -2 ($\alpha < -2$) [11] so that the economically relevant portion of the logit curve, where demand varies meaningfully with price, lies within the feasible domain $p \geq 0$. This prevents the meaningful portion of the demand curve from falling below the clipping point at $p = 0$. Figure 3 illustrates this horizontal shift as $\alpha$ varies.

Figure 4 illustrates the analytical behavior of the logit demand function through its first- and second-order derivatives with respect to price $p$. The first derivative of the demand function ($d'(p)$) remains negative, reflecting the consistent downward slope of the demand curve and confirming that demand decreases as price increases. The steepest point (slope) in this curve corresponds to the region of highest price sensitivity. The first derivative is as follows:

$$d'(p) = -\mu\theta \frac{e^{-(\alpha+\theta p)}}{[1 + e^{-(\alpha+\theta p)}]^2} \tag{2}$$

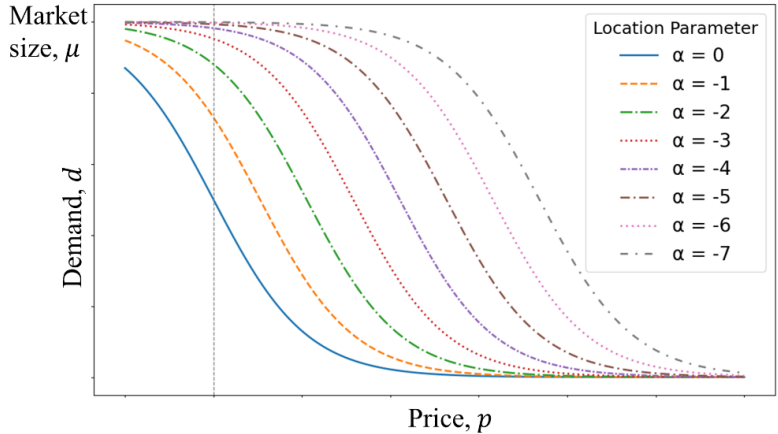

Figure 3: Effect of the location parameter $\alpha$ on the logit demand curve

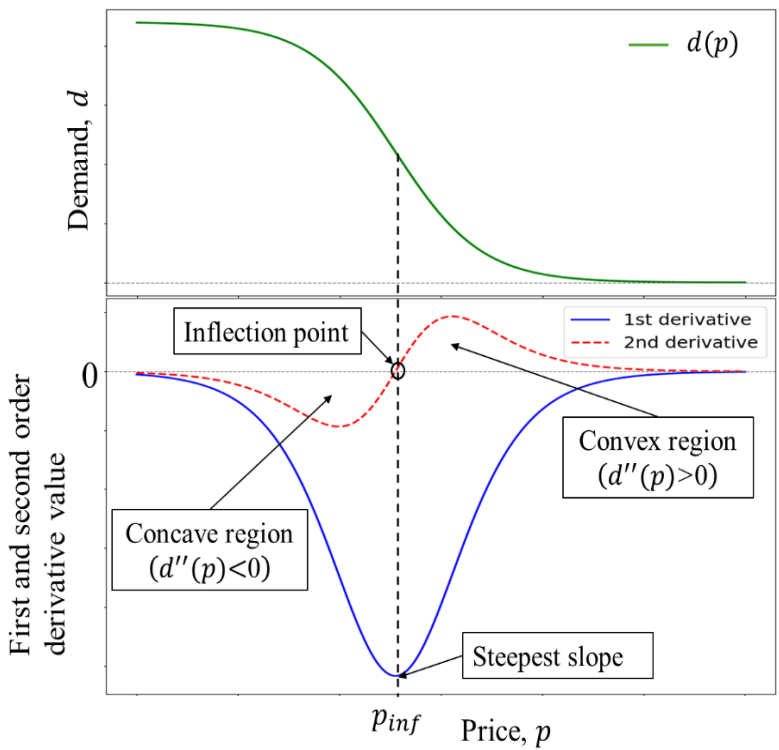

Figure 4: First- and second-order derivatives of the demand function

The second derivative, $d''(p)$, characterizes the curvature of the demand function and helps us understand the function's concavity. This analysis illustrates the regions of the demand curve where consumer response shifts in response to changes in price. The second derivative is as follows:

$$d''(p) = \mu\theta^2 \frac{e^{-(\alpha+\theta p)}[1 - e^{-(\alpha+\theta p)}]}{[1 + e^{-(\alpha+\theta p)}]^3} \tag{3}$$

When the second derivative is negative $(d''(p) < 0)$, the demand curve is concave and when the second derivative is positive $(d''(p) > 0)$, the demand function is convex. The graph clearly shows the concave region and convex region separated by the inflection point where the second derivative becomes zero $(d''(p) = 0)$. Let the price corresponding to the inflection point be denoted by $p_{inf}$. This inflection point is particularly significant, as it corresponds to the highest slope on the demand curve and lies within the elastic region. It is often treated as a reference price and commonly referred to as the market price [9, 11], but it does not generally coincide with the unit elastic point (see section 3.1.1). Further analysis of the inflection point is presented in section 3.1.2.

### 3.1.1 Elasticity of demand

Having characterized the curvature of $d(p)$ through its first and second derivatives in Section 3.1, we now quantify how responsive the demand is to price changes by analyzing the elasticity of the demand function. The sign-aware price elasticity of demand measures the percentage change in demand resulting from a 1% change in price [9, 69], and it is defined as follows:

$$\epsilon(p) = \frac{pd'(p)}{d(p)} \tag{4}$$

For logit demand, let $d(p) = \mu\left(\frac{e^{-(\alpha+\theta p)}}{1+e^{-(\alpha+\theta p)}}\right) = \mu s(p)$, where $s(p) = \frac{1}{1+e^{(\alpha+\theta p)}} \in (0,1)$ denotes the logistic share function. With $\mu > 0$ and $\theta > 0$, this elasticity admits a closed form, derived as follows, using $d'(p) = -\mu\theta \frac{e^{-(\alpha+\theta p)}}{[1+e^{-(\alpha+\theta p)}]^2} = -\mu\theta s(p)[1 - s(p)]$, we obtain the closed form formula for the elasticity of the logit function as follows:

$$\epsilon(p) = \frac{pd'(p)}{d(p)} = -p\theta[1 - s(p)] = -p\theta\left(\frac{e^{(\alpha+\theta p)}}{1 + e^{(\alpha+\theta p)}}\right) \tag{5}$$

The value of $\epsilon(p)$ determines the revenue implications of price changes, a relationship that is standard in microeconomic theory [69, 70]. That is:
- If $-1 < \epsilon(p) < 0$, demand is inelastic, a price increase raises revenue while a price decrease reduces it.
- If $\epsilon(p) = -1$, demand is unit elastic, revenue is maximized at this price.

- If $\epsilon(p) < -1$, demand is elastic, a price increase reduces revenue while a price decrease increases it.

Beyond this general classification, the logit elasticity exhibits further structural properties across the price range. We state these formally below in Proposition 1.

**Proposition 1.** Let $d(p) = \mu s(p)$ where $s(p) = \frac{1}{1+e^{(\alpha+\theta p)}}$ with $\mu > 0$ and $\theta > 0$. For any $p > 0$, the price elasticity of demand $\epsilon(p) = -p\theta[1 - s(p)]$ satisfies $-p\theta < \epsilon(p) < 0$, and is strictly decreasing in $p$.

**Proof.** Given that $e^{(\alpha+\theta p)} > 0$, we have $0 < s(p) < 1$ and $0 < 1 - s(p) < 1$. Thus, for $p > 0$ and $\theta > 0$, we get $-p\theta < \epsilon(p) = -p\theta[1 - s(p)] < 0$. Now, differentiating $\epsilon(p) = -p\theta[1 - s(p)]$ with respect to $p$ we get $\epsilon'(p) = -\theta[1 - s(p)] - p\theta^2 s(p)[1 - s(p)] < 0$. Therefore, $\epsilon(p)$ is strictly decreasing in $p$. ∎

*3.1.2 Inflection point price analysis*

The inflection point is the midpoint of the demand curve, and the price corresponding to that point is commonly referred to as the market price [9, 11]. The inflection point price, $p_{inf}$, is determined by examining the second derivative of the demand function with respect to price $p$. The second derivative of the demand function is as follows:

$$d''(p) = \mu\theta^2 \frac{e^{-(\alpha+\theta p)}[1 - e^{-(\alpha+\theta p)}]}{[1 + e^{-(\alpha+\theta p)}]^3} \tag{6}$$

The inflection point occurs when the second derivative equals zero:

$$d''(p) = 0$$

Because $\mu\theta^2 > 0$ (assuming $\mu, \theta > 0$) and the denominator is strictly positive, for the second derivative to be zero, the numerator of the fraction must equal zero:

$$e^{-(\alpha+\theta p)}[1 - e^{-(\alpha+\theta p)}] = 0 \tag{7}$$

Equation (7) implies two possibilities:

1. $e^{-(\alpha+\theta p)} = 0$, which is not possible because $e^x > 0$ for all real $x$ values.
2. $1 - e^{-(\alpha+\theta p)} = 0$.

Therefore, at the inflection point $1 - e^{-(\alpha+\theta p)} = 0 \Rightarrow e^{-(\alpha+\theta p)} = 1$. Taking the natural logarithm of both sides and simplifying the equation, we get the following formula for the inflection point price $p_{inf}$:

$$p_{inf} = -\frac{\alpha}{\theta} \tag{8}$$

While inflection point of the logit demand function has historically been used as a reference for strategic pricing, we show in Proposition 2 that this point does not necessarily correspond to the price that maximizes revenue. Rather, it corresponds to the point where the slope of the demand curve is steepest, i.e., demand sensitivity to price is highest.

**Proposition 2.** Let $d(p) = \mu s(p)$ where $s(p) = \frac{1}{1+e^{(\alpha+\theta p)}}$ with $\mu > 0$, $\alpha < -2$, and $\theta > 0$. Then the price elasticity of demand, $\epsilon(p) = -p\theta[1 - s(p)]$, at the inflection point price $p_{inf}$ does not satisfy $\epsilon(p_{inf}) = -1$, and hence $p_{inf}$ is not the revenue-maximizing price.

**Proof.** From the second derivative of the demand function $d(p)$, we have $e^{-(\alpha+\theta p_{inf})} = 1$, therefore, $s(p_{inf}) = \frac{1}{1+e^{(\alpha+\theta p)}} = \frac{1}{2}$. At the inflection point $p_{inf} = -\frac{\alpha}{\theta}$ then $\epsilon(p_{inf}) = -p_{inf}\theta[1 - s(p_{inf})] = -\left(-\frac{\alpha}{\theta} \times \theta\right)\frac{1}{2} = \frac{\alpha}{2}$. This equals $-1$ only when $\alpha = -2$, given that $\alpha < -2$ then the inflection point cannot coincide with the unit-elastic point and $p_{inf}$ is not the revenue maximizing price. ∎

Furthermore, the elasticity level at the inflection point is based solely on the value of $\alpha$. Only when $\alpha = -2$, the price elasticity, $\epsilon(p_{inf})$, becomes one and the inflection point price coincides with the revenue-maximizing price. However, for most practical applications $\alpha < -2$, larger values (closer to zero) shift the demand curve to the left, meaning that a significant portion of the market cannot be captured unless we offer negative prices, as can be seen in Figure 3, which is not realistic.

### 3.2 Exact pricing algorithm

In this subsection, we develop an Exact Pricing Algorithm (EPA) that identifies the revenue-maximizing price $p^*$ by solving the first-order condition of the revenue function. The revenue, denoted by $R$, for a product is expressed as:

$$R = pd(p) \tag{9}$$

Next, we derive the revenue function with respect to $p$:

$$R'(p) = \frac{d}{dp}(pd(p)) = d(p) + pd'(p) \tag{10}$$

Substituting the logistic demand function into equation (10) yields:

$$R'(p) = \mu \left[ \frac{e^{-(\alpha+\theta p)}}{1 + e^{-(\alpha+\theta p)}} - p\theta \frac{e^{-(\alpha+\theta p)}}{(1 + e^{-(\alpha+\theta p)})^2} \right] \tag{11}$$

Setting the derivative of the revenue function to zero leads to the critical point for revenue maximization allowing us to obtain the equation for the optimal price $p^*$ as:

$$p^*\theta - 1 = e^{-(\alpha+\theta p^*)} \Rightarrow p^* = \frac{e^{-(\alpha+\theta p^*)} + 1}{\theta} \tag{12}$$

Multiplying both sides of the equation (12) by $e^{(\alpha+\theta p)}$, we have the following equation:

$$(p^*\theta - 1)e^{(\alpha+\theta p^*)} = 1 \tag{13}$$

Equation (13) is transcendental because $p^*$ is both inside and outside the exponential function, therefore, it cannot be solved by algebraic techniques. To address this complexity, we utilize the Lambert $W$ function, also known as the product logarithm, which is a special function used to solve equations where a variable appears both inside and outside an exponential term, such as in the form $ze^z = x \Rightarrow z = W(x)$. For more details on the Lamber $W$ function, see Corless et al. [71].

By manipulating the derived equation for $p^*$, we rewrite it into a form compatible with the Lambert $W$ function: $e^{\alpha+\theta p^*}(\theta p^* - 1) = 1$. Letting $z = \theta p^* - 1 \Rightarrow \theta p^* = z + 1 \Rightarrow p^* = \frac{z+1}{\theta}$. Now, the Lambert $W$ function form of equation (13) is $ze^z = e^{-(\alpha+1)} \Rightarrow z = W(e^{-(\alpha+1)})$. Substituting back $z = \theta p^* - 1$, we solve for the optimal price $p^*$ as:

$$\theta p^* - 1 = W(e^{-(\alpha+1)})$$
$$\Rightarrow p^* = \frac{1}{\theta}\left[1 + W(e^{-(\alpha+1)})\right] \tag{14}$$

This representation avoids ad-hoc iterative root-finding for $p^*$ and is numerically stable to evaluate using standard implementations of $W$. After determining the optimal price $p^*$ using equation (14), the final revenue function can be expressed by substituting $p^*$ into the revenue equation (9). Now, the revenue $R$ is calculated as:

$$R(p^*) = p^* d(p^*)$$

$$\Rightarrow R(p^*) = \left[p^* \times \left[\mu\left(\frac{e^{-(\alpha+\theta p^*)}}{1 + e^{-(\alpha+\theta p^*)}}\right)\right]\right] \tag{15}$$

**Proposition 3.** The optimal price obtained by EPA, $p^*$, satisfies $\epsilon(p^*) = -1$.

**Proof.** From the first-order condition of the revenue function, we have:
$$p^*\theta - 1 = e^{-(\alpha+\theta p^*)} \Rightarrow p^*\theta = 1 + e^{-(\alpha+\theta p^*)}$$
Now, substitute $p^*\theta = 1 + e^{-(\alpha+\theta p^*)}$ into equation (5):
$$\epsilon(p^*) = -(1 + e^{-(\alpha+\theta p^*)})\left(\frac{e^{(\alpha+\theta p^*)}}{1+e^{(\alpha+\theta p^*)}}\right),$$

$$\epsilon(p^*) = -\frac{1+e^{-(\alpha+\theta p)}}{1+e^{-(\alpha+\theta p)}},$$
$$\epsilon(p^*) = -1.$$

Hence, $p^*$ is the optimal price. ∎

Now to understand the broader behavior of the revenue function $R(p)$ over the entire feasible price domain, we must verify that the critical point $p^*$ identified by EPA is not only a stationary point but in fact the unique global maximizer of revenue. The following theorem formalizes these properties by characterizing the monotonicity and unimodality of $R(p)$, thereby guaranteeing both existence and uniqueness of the optimal price $p^*$.

**Theorem 1.** Let $d(p)$ be a logit demand function with $\mu > 0$, $\theta > 0$, $\alpha < -2$, and let $R = pd(p)$ be the corresponding revenue function. There exists an optimal price $p^*$ such that:
- $R(p)$ attains its global maximum at $p^*$;
- $R'(p) > 0$ for all $p < p^*$, and $R'(p) < 0$ for all $p > p^*$;
- Consequently, $R(p)$ is strictly increasing on $(0, p^*)$, strictly decreasing on $(p^*, \infty)$, and unimodal.

**Proof.** To determine the optimal price $p^*$ that maximizes the revenue $R(p) = pd(p)$, we identify the critical point by setting the first derivative of the revenue function to zero:

$$R'(p^*) = d(p^*) + p^* d'(p^*) = 0$$

Solving the derivative of the revenue function $R'(p^*)$ leads to a transcendental equation, $p^*\theta - 1 = e^{-(\alpha+\theta p^*)}$. This transcendental equation is solved using the Lambert $W$ function to find a critical point $p^*$ where $p^* = \frac{1}{\theta}\left[1 + W\left(e^{-(\alpha+1)}\right)\right]$.

Next, we examine uniqueness and monotonicity. Price elasticity of demand is $\epsilon(p) = \frac{pd'(p)}{d(p)}$. Then substituting $\epsilon(p)$ into the derivative of the revenue function $R'(p)$ can be expressed as $R'(p) = d(p) + pd'(p) = d(p)(1 + \epsilon(p))$ where $\epsilon(p) = -p\theta[1 - s(p)]$. From Proposition 1, $\epsilon(p)$ is strictly decreasing in $p$. From Proposition 3, the optimal solution satisfies $\epsilon(p^*) = -1$. Elasticity decreases monotonically and elasticity at optimal price $\epsilon(p^*) = -1$, it follows that $\epsilon(p) > -1$ for all $p < p^*$, so $1 + \epsilon(p) > 0$ yielding $R'(p) > 0$. Now, for all $p > p^*$, and $\epsilon(p) < -1$, so $1 + \epsilon(p) < 0$ yielding $R'(p) < 0$. Thus, $R(p)$ is strictly increasing on $(0, p^*)$, and strictly decreasing on $(p^*, \infty)$.

Moreover, as $p \to \infty$, $s(p) \to 0$, so $d(p) \to 0$ and $\epsilon(p) \to -\infty$. Hence $R'(p) = d(p)(1 + \epsilon(p)) \to 0^-$. This establishes that the slope of the revenue function $R'(p)$ approaches zero but stays negative as $p \to \infty$. Since revenue is strictly increasing before $p^*$, strictly decreasing after $p^*$, the function has exactly one global maximum at $p^*$. This establishes that $R(p)$ is unimodal with unique maximizer $p^*$. ∎

Figure 5 shows the behavior of $R(p)$, $R'(p)$, and $R''(p)$ confirming Theorem 1 where the revenue function $R(p)$ peaks at $p^*$. Having established in Theorem 1 that $p^*$ is the unique global maximizer of revenue and that $R(p)$ decreases monotonically beyond this point, we now compare between the revenue obtained from the optimal price $p^*$ and the inflection-point price $p_{inf}$. In Proposition 4, we will prove that $p_{inf} > p^*$.

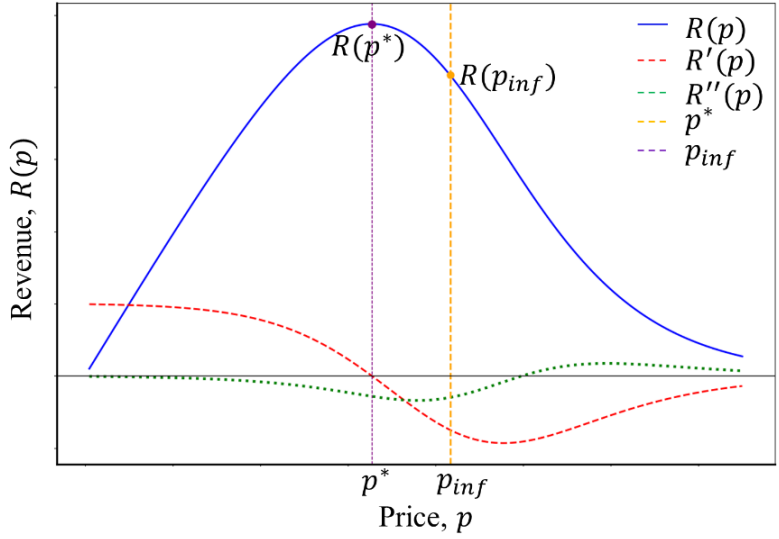

Figure 5: Revenue function and its second derivative

**Proposition 4.** Let $d(p)$ be a logit demand function with $\mu > 0$, $\theta > 0$, $\alpha < -2$, and let $R = pd(p)$ be the corresponding revenue function. Then the revenue-maximizing price $p^*$ satisfies $p^* < p_{inf}$, where $p_{inf}$ is the price associated with the inflection point of $d(p)$.

**Proof.** The inflection point of the logit demand function occurs at $1 - e^{-(\alpha+\theta p)} = 0$, thus the inflection point price satisfies $\alpha + \theta p_{inf} = 0$ which implies $p_{inf} = -\frac{\alpha}{\theta}$. To analyze the revenue behavior around this point, we consider the first derivative of the revenue function. The derivative of revenue function with respect to price $p$ is $R'(p) = d(p) + pd'(p)$. Now, let $u = \alpha + \theta p_{inf}$, which simplifies to $u = \alpha + \theta\left(-\frac{\alpha}{\theta}\right) = 0$. This gives $e^{-u} = e^{-0} = 1$. Substituting into the demand function, we obtain $d(p_{inf}) = \mu \frac{1}{1+1} = \frac{\mu}{2}$. The derivative of the demand function is $d'(p) = -\mu\theta \frac{e^{-(\alpha+\theta p)}}{[1+e^{-(\alpha+\theta p)}]^2}$, and at the inflection point where $e^{-(\alpha+\theta p)} = 1$, this simplifies to $d'(p_{inf}) = -\mu\theta \frac{1}{[1+1]^2} = -\frac{\mu\theta}{4}$. Now, using these values we compute the derivative of the revenue function at $p_{inf}$ as follows:

$$R'(p_{inf}) = d(p_{inf}) + p_{inf}d'(p_{inf}) = \frac{\mu}{2} + \left(-\frac{\alpha}{\theta}\right)\left(-\frac{\mu\theta}{4}\right) = \mu\left(\frac{1}{2} + \frac{\alpha}{4}\right)$$

Under the assumptions that $\alpha < -2$ and $\mu > 0$, it follows that $R'(p_{inf}) < 0$. This means the slope of the revenue function is negative at the inflection point. By Theorem 1 the revenue strictly increasing from $(0, p^*)$ and strictly decreasing $(p^*, \infty)$. Therefore, $p_{inf} \in (p^*, \infty)$ and $p^* < p_{inf}$.
∎

Proposition 4 has an interesting implication: firms selling products that experience price-sensitive demand, following a logit function, can reduce their prices below the commonly used inflection point price to maximize revenue and, at the same time, gain a competitive advantage. This outcome aligns with customer welfare, as the revenue-maximizing strategy coincides with more affordable prices.

To implement the proposed methodology for revenue maximization, we present the following step-by-step algorithm for computing the optimal price.

---
**Exact Pricing Algorithm**

**Step 1**: Initialize parameters for revenue maximization
- Identify function parameters including location parameter ($\alpha$), market size ($\mu$), and price sensitivity ($\theta$).
- Define the demand and revenue functions, $d(p) = \mu \left( \frac{e^{-(\alpha+\theta p)}}{1+e^{-(\alpha+\theta p)}} \right)$, $R(p) = pd(p)$.

**Step 2**: Solve for optimal price $p^*$

**Sub-step 2.1**:
- Derive the revenue function and set the derivative equal to zero: $R'(p) = d(p) + pd'(p) = 0$.
- Substituting $d(p)$, the first-order condition simplifies to: $(p\theta - 1)e^{(\alpha+\theta p)} = 1$.

**Sub-step 2.2**: Apply the Lambert $W$ function:
- Let $z = p\theta - 1$, then: $ze^z = e^{-(\alpha+1)} \Rightarrow z = W(e^{-(\alpha+1)})$.
- Solve for optimal price: $p^* = \frac{1}{\theta}[1 + W(e^{-(\alpha+1)})]$.

**Sub 3**: Compute the demand and revenue at $p^*$
- Compute the demand $d(p^*)$ at the optimal price $p^*$ using $d(p^*) = \mu \left( \frac{e^{-(\alpha+\theta p^*)}}{1+e^{-(\alpha+\theta p^*)}} \right)$.
- Compute the revenue $R(p^*)$ at the optimal price $p^*$ using $R = p^* d(p^*)$.

---

## 3.3 Revenue and price ratios

To further quantify the improvement obtained by EPA over the commonly used inflection point pricing $p_{inf}$, let $R^* = R(p^*)$ and $R_{inf} = R(p_{inf})$. We define the revenue ratio in equation (16) as follows:

$$\frac{R^*}{R_{inf}} = \frac{1 + W(e^{-(\alpha+1)})}{-\alpha} \times e^{-\theta(p^* - p_{inf})} \frac{(1 + e^{-(\alpha+\theta p_{inf})})}{(1 + e^{-(\alpha+\theta p^*)})} \quad (16)$$

This expression measures the relative revenue achieved by pricing at the optimal point $p^*$ versus at the inflection point $p_{inf}$. Note that, by definition, $\frac{R^*}{R_{inf}} \geq 1$, reflecting the superiority of the optimal price over the inflection point benchmark. In addition to the revenue ratio, we can define the ratio of the optimal price $p^*$ to the inflection point price $p_{inf}$ in equation (17) as follows:

$$\frac{p^*}{p_{inf}} = \frac{1 + W(e^{-(\alpha+1)})}{-\alpha} \tag{17}$$

This ratio depends on both the Lambert $W$ function and the location parameter $\alpha$. Since $W(e^{-(\alpha+1)}) < 1$ when $\alpha < -2$, we obtain $\frac{p^*}{p_{inf}} < 1$, confirming that the inflection-point price consistently overestimates the true optimal price. Therefore, the revenue ratio and the price ratio further support the adoption of our exact pricing algorithm.

## 4. Numerical analysis

To evaluate the practical implications of the proposed EPA under the logit demand function, this section presents a set of numerical experiments. The objective is to validate the analytical results and compare the revenue-maximizing price $p^*$ with the commonly used inflection point price $p_{inf}$. We illustrate the performance of the closed-form solution through visualizations of the revenue and demand functions and conduct a sensitivity analysis to examine how different parameter settings affect demand behavior, optimal pricing decisions, and resulting revenue. We begin our analysis using a representative baseline set of parameters as follows:

- Market size $(\mu) = 1,000$
- Location parameter $(\alpha) = -6$
- Price sensitivity parameter $(\theta) = 0.30$

### 4.1 Optimal price computation by EPA

Using the parameter values above, the optimal price $p^*$ is computed by the proposed EPA:

$$p^* = \frac{1}{\theta}[1 + W(e^{-(\alpha+1)})]$$

Substituting $\alpha = -6, \theta = 0.30$, the optimal price $p^*$ is calculated as:

$$p^* = \frac{1}{0.30}[1 + W(e^{-(-6+1)})]$$

Evaluating the Lambert $W$ function numerically yields:

$$W(e^5) \approx 3.69344 \Rightarrow p^* \approx 15.64$$

Next, the demand at the optimal price is:

$$d(p^*) = \mu\left(\frac{e^{-(\alpha+\theta p^*)}}{1 + e^{-(\alpha+\theta p^*)}}\right) = 1{,}000 \times \frac{e^{-(-6+0.30\times 15.64)}}{1 + e^{-(-6+0.30\times 15.64)}} \approx 786.94$$

The corresponding revenue at the optimal price is:

$$R(p^*) = p^* d(p^*) \approx \$15.64 \times 786.94 \approx \$12{,}311.47$$

### 4.2 Inflection point price

Mathematically, the inflection point occurs when the exponent in the denominator of the logit function equals zero, presented in equation (8) as $p_{inf} = -\frac{\alpha}{\theta}$. Substituting $\alpha = -6$ and $\theta = 0.30$ gives $p_{inf} = -\frac{-6}{0.30} = \$20.00$.

At this point, the expected demand equals half of the market size, $d(p_{\text{inf}}) = \frac{\mu}{2} = 500$, and revenue is

$$R(p_{inf}) = p_{inf} \frac{\mu}{2} = \$20.00 \times 500 = \$10{,}000.$$

Figure 6 and Figure 7 illustrate the demand and revenue behavior over the price range. As implied by the logit structure, demand decreases monotonically with price, with the inflection point marking the steepest decline.

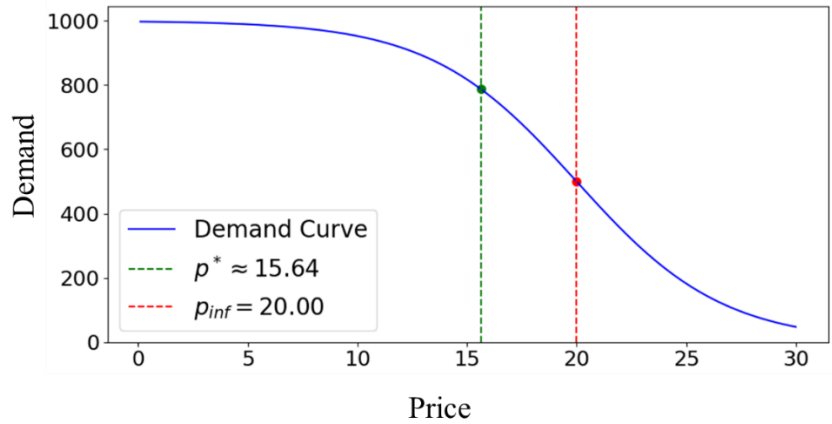

Figure 6: Demand vs Price curve

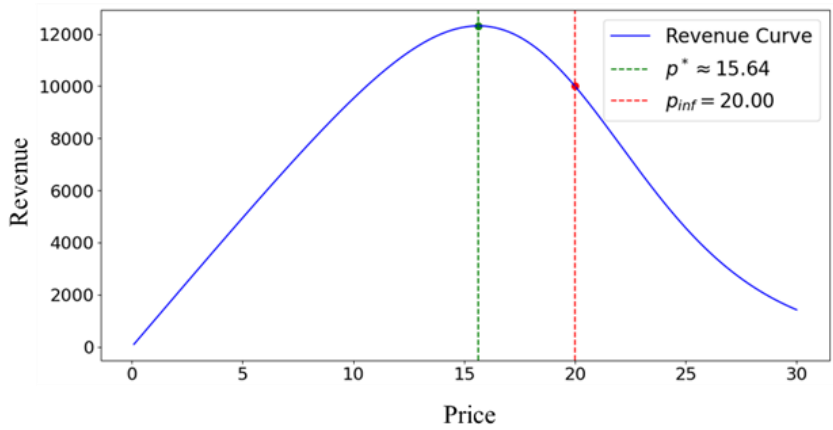

Figure 7: Revenue vs price curve

These visualizations confirm that although the inflection point indicates where consumer sensitivity to price is highest, the optimal price must balance that sensitivity with the actual volume of sales, which is precisely what the proposed closed-form EPA solution achieves. The revenue curve initially rises as both price and demand contribute positively, reaches its maximum at the optimal price $p^*$, and then falls as demand diminishes. Notably, $p^*$ lies to the left of the inflection point price $p_{inf}$. At $p_{inf}$, demand is already decaying rapidly, producing lower revenue than at the optimum.

### 4.3 Elasticity of demand

Figure 8 illustrates the behavior of the price elasticity of demand, $\epsilon(p)$, over a range of prices, computed using the closed-form expression given in equation (5):

$$\epsilon(p) = \frac{pd'(p)}{d(p)} = -p\theta[1 - s(p)] = -p\theta\left(\frac{e^{\alpha+\theta p}}{1 + e^{\alpha+\theta p}}\right)$$

In the plotted example, we use $\alpha = -6$ and $\theta = 0.30$. Two benchmark prices are highlighted:

1. Revenue maximizing price $p^*$ is the point at which elasticity equals unity, i.e., $\epsilon(p^*) = -1$. This condition follows from the first-order condition of revenue maximization:
$$R'(p) = d(p) + pd'(p) = 0 \Rightarrow \epsilon(p) = -1.$$
Using the closed-form expression from Section 4.1, we obtain $p^* \approx 15.64$ for the current parameterization, which is confirmed visually in Figure 8 by the intersection of the elasticity curve with the unit-elastic threshold.

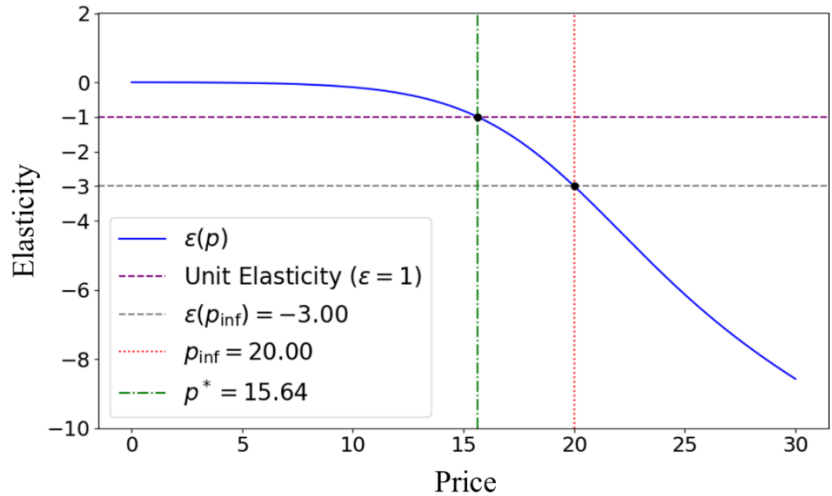

Figure 8: Price elasticity of demand across prices

2. Inflection point price $p_{inf}$, derived from the second derivative $d''(p) = 0$, occurs at $p_{inf} = 20$ and the corresponding elasticity is:
$$\epsilon(p_{inf}) = -3,$$
which is consistent with Proposition 2. This value is shown by the dashed line, clearly distinct from the unit-elastic reference.

Figure 8 confirms that $p_{inf} \neq p^*$, reinforcing the analytical results in Proposition 2. Specifically, while $p^*$ identifies the point of maximum revenue, $p_{inf}$ only reflects a geometric property of the demand function. Additionally, Figure 8 illustrates that the elasticity decreases monotonically with price, consistent with Proposition 1. This monotonicity ensures that the unit-elastic point is unique and occurs before the inflection point when $\alpha < -2$, as is the case here $\alpha = -6$.

### 4.4 Sensitivity analysis

We conduct a sensitivity analysis to understand the influence of parameter values on the optimal pricing decision and the corresponding revenue outcomes. The parameter ranges are selected to reflect realistic and interpretable behavior of the logit demand function. The interval $\alpha \in [-7, -3]$ spans high to moderate baseline demand while remaining within the empirically relevant region $\alpha < -2$. The range $\theta \in [0.1, 0.5]$ represents low to moderately high price sensitivity and corresponds to demand environments where consumers respond meaningfully to price changes. Together, these ranges allow the sensitivity analysis to cover a broad spectrum of plausible market conditions.

For each scenario, we calculate the optimal price using the proposed EPA and compute the associated demand and revenue. This analysis allows us to observe trends and interpret pricing behaviors. Table 2 summarizes the outcomes of the sensitivity analysis across various combinations of $\alpha$ and $\theta$ values. Each row presents the computed optimal price $p^*$, the demand at the optimal price $d(p^*)$, the resulting revenue $R(p^*)$, the inflection point price $p_{inf}$, the demand at the inflection point price $d(p_{inf})$, and the revenue at the inflection point $R(p_{inf})$. From the table, several insights emerge:

- As $\alpha$ increases (from –7 to –3), both $p^*$ and $p_{inf}$ decrease. This is expected, as larger $\alpha$ shifts the demand curve to the left, leading to lower prices.
- As $\theta$ increases, indicating higher price sensitivity, both $p^*$ and $p_{inf}$ decrease. This is due to the accelerating decline in demand with increasing prices.
- For all scenarios, $p^* < p_{inf}$ and $R(p^*) > R(p_{inf})$. This reinforces the core contribution of this study: maximum revenue is not attained at the inflection point but rather at a lower price where the trade-off between price and demand is optimized. It also supports Theorem 1 and Proposition 4.

Table 2 shows that the demand corresponding to the inflection point price $p_{inf}$ remains constant at $\mu/2$ (i.e., 500 units when $\mu = 1000$) across all parameter settings. This invariance follows directly from Proposition 4, which establishes that at the inflection point the logit share equals one-half, so that $d(p_{inf}) = \frac{\mu}{2}$. In contrast, the optimal price $p^*$ is determined by maximizing revenue, which depends on the specific trade-off between price and quantity demanded. The optimal price balances higher unit price against reduced demand to achieve the maximum possible revenue. Consequently, demand at optimal price $d(p^*)$ changes with variations in $\alpha$ and $\theta$, reflecting how these parameters shape consumers' willingness to pay.

Table 2: Summary table for various combinations of $\alpha$ and $\theta$

| Case | $\alpha$ | $\theta$ | $\mu$ | $p^*$ | $d(p^*)$ | $R(p^*)$ | $p_{inf}$ | $d(p_{inf})$ | $R(p_{inf})$ |
|---|---|---|---|---|---|---|---|---|---|
| 1 | -7 | 0.1 | 1,000 | $54.97 | 818.1 | $44,966.6 | $70 | 500 | $35,000 |
| 2 | -5 | 0.1 | 1,000 | $39.26 | 745.3 | $29,262.7 | $50 | 500 | $25,000 |
| 3 | -4 | 0.1 | 1,000 | $32.08 | 688.3 | $22,079.4 | $40 | 500 | $20,000 |
| 4 | -3 | 0.1 | 1,000 | $25.57 | 608.9 | $15,571.5 | $30 | 500 | $15,000 |
| 5 | -7 | 0.3 | 1,000 | $18.32 | 818.1 | $14,988.9 | $23.33 | 500 | $11,666.70 |
| 6 | -5 | 0.3 | 1,000 | $13.09 | 745.3 | $9,754.2 | $16.67 | 500 | $8,333.30 |
| 7 | -4 | 0.3 | 1,000 | $10.69 | 688.3 | $7,359.8 | $13.33 | 500 | $6,666.70 |
| 8 | -3 | 0.3 | 1,000 | $8.52 | 608.9 | $5,190.5 | $10 | 500 | $5,000 |
| 9 | -7 | 0.5 | 1,000 | $10.99 | 818.1 | $8,993.3 | $14 | 500 | $7,000 |
| 10 | -5 | 0.5 | 1,000 | $7.85 | 745.3 | $5,852.5 | $10 | 500 | $5,000 |
| 11 | -4 | 0.5 | 1,000 | $6.42 | 688.3 | $4,415.9 | $8 | 500 | $4,000 |
| 12 | -3 | 0.5 | 1,000 | $5.11 | 608.9 | $3,114.3 | $6 | 500 | $3,000 |

While Table 2 provides detailed numerical results for different combinations of $\alpha$ and $\theta$, visualizations can further clarify the underlying patterns. To this end, Figure 9 and Figure 10 present graphical summaries based on the data presented in Table 2. Figure 9 illustrates how the optimal price $p^*$ and revenue at the optimal price $R(p^*)$ vary with the price sensitivity parameter $\theta$, for different values of the location parameter $\alpha$. Each line represents a distinct $\alpha$ level. The key observations are as follows:

- As the price sensitivity parameter $\theta$ increases, the optimal price $p^*$ decreases monotonically for all $\alpha$ values.
- The revenue at the optimal price $R(p^*)$, while also decreasing $\theta$, does so at a slightly slower rate, reflecting the trade-off between lower prices and higher quantities sold.
- The location parameter $\alpha$ exerts a strong influence on both pricing and revenue: more negative $\alpha$ values yield higher $p^*$ and greater revenue potential. However, the gap between different $\alpha$ curves shrink as $\theta$ increases.

Finally, Figure 9 demonstrates the trade-off between demand sensitivity and price, highlighting that aggressive pricing (low $p^*$) is essential when consumers are highly responsive to price (large $\theta$).

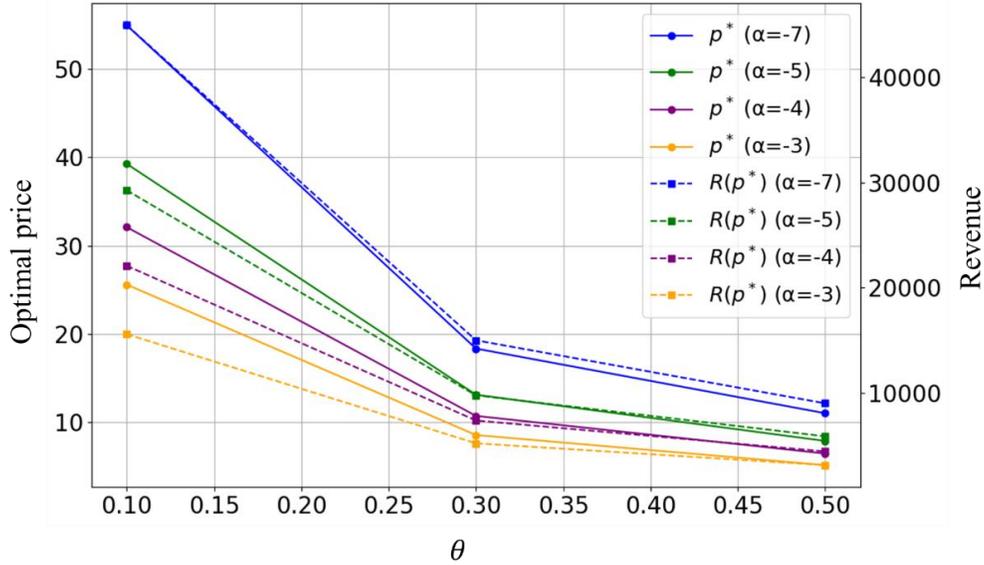

Figure 9: Effect of $\theta$ on optimal price and revenue for various $\alpha$ values

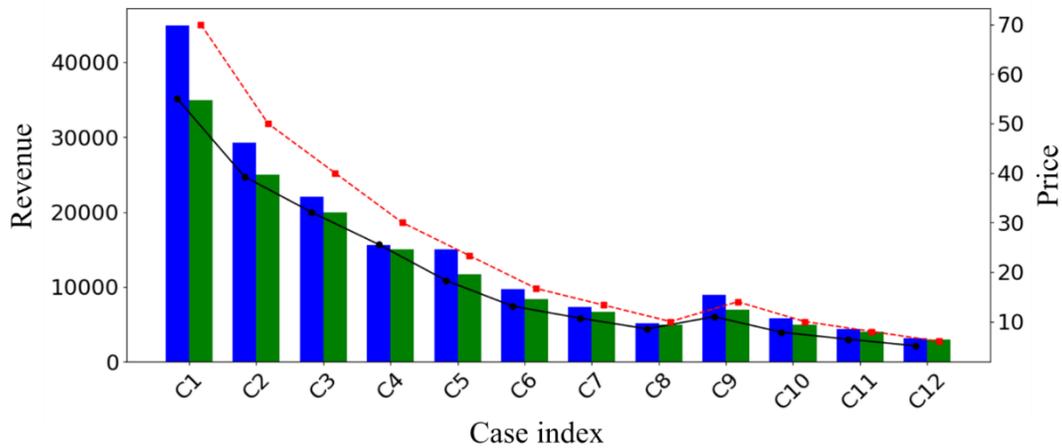

Figure 10: Comparison of optimal and inflection prices with corresponding revenues

Figure 10 compares $p^*$ and $R(p^*)$ against $p_{inf}$ and $R(p_{inf})$ for each parameter setting using values from Table 2. Note that Cases 1-12 are not ordered or sequential; they simply index the combinations of $\alpha$ and $\theta$ used in Table 2. The key observations are as follows: for every case, the bar for $R(p^*)$ is consistently taller than that for $R(p_{inf})$. The secondary Y-axis highlights that $p^*$ is lower than $p_{inf}$ for all cases. An important trend visible in Figure 10 is that the gap between $p^*$ and $p_{inf}$ gradually decreases as $\alpha$ increases (from $-7$ to $-3$) and as $\theta$ increases. This is consistent with the theoretical results that $p^*$ and $p_{inf}$ coincide when $\alpha = -2$. Similarly, an increase in $\theta$ accelerates the decline of the demand curve, bringing $p^*$ and $p_{inf}$ closer together.

The sensitivity analysis also reveals the distinct influence of the location parameter $\alpha$ on the demand curve and optimal pricing outcomes. Figure 11 further illustrates the sensitivity of the logit demand function to changes in the location parameter $\alpha$. As $\alpha$ becomes more negative, the demand curve shifts progressively to the right, indicating that consumers are willing to accept higher prices before demand begins to decline. This shift reflects a higher baseline valuation of the product and a greater tolerance for price increases. In contrast, less negative $\alpha$ values (e.g., $-3$ or $-4$) result in leftward-shifted curves, showing lower demand at low price points and a steeper drop-off as price increases, consistent with markets characterized by stronger price sensitivity and lower perceived product value.

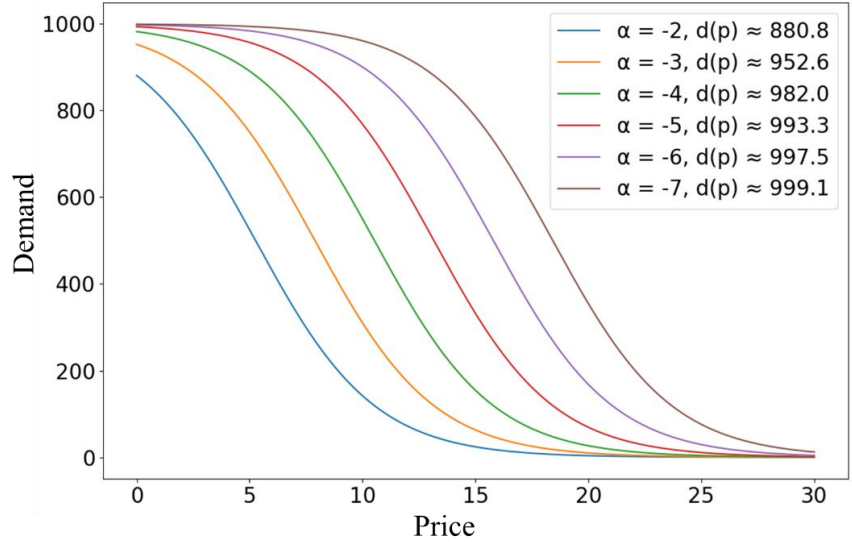

Figure 11: Logit demand function for varying $\alpha$ values

These results complement the numerical findings in Table 2, where more negative $\alpha$ values lead to higher optimal prices and greater achievable demand at the revenue-maximizing price. Together, they demonstrate that $\alpha$ governs not only the horizontal position of the demand curve but also the effective range of viable prices. From a managerial perspective, accurately estimating $\alpha$ is essential: more negative $\alpha$ values represent premium or high-utility products with relatively stable demand across prices, whereas less negative $\alpha$ values describe commodity-like items requiring more competitive pricing strategies.

### 4.5 Analysis of revenue and price ratios

The revenue ratio $\frac{R^*}{R_{inf}}$ and price ratio $\frac{p^*}{p_{inf}}$ provide insights for managers and pricing strategists. These ratios, computed using equations (16) and (17), respectively, are summarized in Table 3. The revenue ratio consistently exceeds 1 across all parameter scenarios, confirming that pricing at the theoretically derived optimal price $p^*$ yields higher revenue than relying on the simpler inflection point price $p_{inf}$, which only reflects maximum price sensitivity. From a practical perspective, adopting the optimal price not only increases revenue but also enhances market competitiveness by attracting more price-sensitive customers, improving inventory turnover, and reducing waste.

Table 3: Evaluation of revenue ratio $\frac{R^*}{R_{inf}}$ and price ratio $\frac{p^*}{p_{inf}}$ across different parameter settings

| Case | $\alpha$ | $\theta$ | $\mu$ | $p^*$ | $p_{inf}$ | $\frac{R^*}{R_{inf}}$ | $\frac{p^*}{p_{inf}}$ |
|---|---|---|---|---|---|---|---|
| 1 | -7 | 0.1 | 1000 | 54.97 | 70 | 1.285 | 0.785 |
| 2 | -5 | 0.1 | 1000 | 39.26 | 50 | 1.171 | 0.785 |
| 3 | -4 | 0.1 | 1000 | 32.08 | 40 | 1.104 | 0.802 |
| 4 | -3 | 0.1 | 1000 | 25.57 | 30 | 1.038 | 0.852 |
| 5 | -7 | 0.3 | 1000 | 18.32 | 23.33 | 1.285 | 0.785 |
| 6 | -5 | 0.3 | 1000 | 13.09 | 16.67 | 1.171 | 0.785 |
| 7 | -4 | 0.3 | 1000 | 10.69 | 13.33 | 1.104 | 0.802 |
| 8 | -3 | 0.3 | 1000 | 8.52 | 10 | 1.038 | 0.852 |
| 9 | -7 | 0.5 | 1000 | 10.99 | 14 | 1.285 | 0.785 |
| 10 | -5 | 0.5 | 1000 | 7.85 | 10 | 1.171 | 0.785 |
| 11 | -4 | 0.5 | 1000 | 6.42 | 8 | 1.104 | 0.802 |
| 12 | -3 | 0.5 | 1000 | 5.11 | 6 | 1.038 | 0.852 |

## 5. Conclusions and managerial insights

This study presents an exact, closed-form, and analytically grounded pricing algorithm that maximizes revenue under the logit demand function, representing a significant advancement over traditional heuristic-based approaches. By deriving the first-order condition of the revenue function and solving the resulting transcendental equation using the Lambert $W$ function, we establish a closed-form expression for the revenue-maximizing price. This formulation eliminates the need for iterative or approximate numerical approaches, enhancing computational efficiency and providing a generalizable pricing strategy.

Theoretical results and numerical experiments both confirm that the commonly referenced inflection point price, where demand sensitivity is maximized, does not coincide with the revenue-maximizing price. Instead, the optimal price consistently lies to the left of the inflection point on the revenue curve, where the balance between price and demand volume yields the highest revenue at a lower selling price. It is worth noting that, due to the unimodal nature of the revenue function, the same revenue obtained at the inflection point price can also be achieved at a lower price, even lower than the optimal price, allowing firms to gain a competitive advantage in the market while maintaining profitability, however, at a lower revenue point than the optimal revenue guaranteed by the optimal price.

Across all examined parameter scenarios, the revenue ratio $\frac{R^*}{R_{inf}}$ remained greater than one, reinforcing the superiority of the EPA-derived price. On average, the EPA approach achieved 15% higher revenue while setting prices about 20% lower than those obtained from the inflection-point heuristic. This aligns with the observation that the price ratio $\frac{p^*}{p_{inf}}$ was always less than one, confirming that the inflection point price tends to overestimate the optimal price, potentially resulting in lost sales and reduced consumer welfare. Sensitivity analysis revealed that as both the location and the price-sensitivity parameters, $\alpha$ and $\theta$ respectively, increase, the gap between the optimal price $p^*$ and the inflection point price $p_{inf}$ narrows, illustrating the convergence of optimal and heuristic pricing in less elastic markets. Additionally, as $\alpha$ increases the improvement in revenue becomes smaller, i.e. $\left(\frac{R^*}{R_{inf}}\right)$ gets closer to 1. Similarly, the ratio between optimal and inflection point prices become higher, i.e. $\left(\frac{p^*}{p_{inf}}\right)$ gets closer to 1, but at a slower rate. Therefore, for cheaper products, the benefit to the customer from optimizing the price exceeds the benefit to the retailer.

From a managerial perspective, the proposed EPA offers a robust, practical, and data-driven tool that can be directly applied without requiring extensive experimentation or computational effort. It enables firms to set prices that not only align with realistic consumer behavior but also maximize profitability while enhancing customer value. By incorporating revenue and price ratios as benchmark indicators, managers can assess their current pricing strategies and make informed adjustments in response to varying market conditions. This analytical transparency supports faster decision-making and enhances their long-term financial performance.

Future research can extend this work in several directions. First, the model could be generalized to multi-product settings to capture substitution cannibalization effects. Second, incorporating dynamic consumer behavior learning mechanisms would allow the algorithm to adapt as new demand data become available. Third, integrating EPA with inventory and supply chain decision models could enable holistic revenue optimization across operational stages. Finally, exploring an age-and-price-dependent logit demand function, which explicitly accounts for product aging effects alongside price sensitivity, represents a promising avenue for further enhancing the practical relevance and managerial value of the proposed pricing approach, particularly for perishable and time-sensitive goods.